\documentclass[10pt]{amsart}
\title[]{
Closed planar curves without inflections}

\date{March 17, 2011.}
\usepackage[dvips]{graphicx} 

\usepackage[mathscr]{eucal}
\usepackage{mathrsfs}
\usepackage{amscd}
\usepackage{amsfonts}
\usepackage{amsmath,amsthm,amssymb}
\usepackage{latexsym}
\usepackage{bm}
\usepackage[dvips]{graphics}
\numberwithin{equation}{section}

\usepackage{amsthm}
\theoremstyle{plain}
 \newtheorem{theorem}{Theorem}[section]
 \newtheorem*{theorem*}{Theorem}

 \newtheorem*{lemma*}{Lemma}
 \newtheorem{proposition}[theorem]{Proposition}
 
 \newtheorem*{fact*}{Fact}

 \newtheorem{lemma}[theorem]{Lemma}
 \newtheorem{corollary}[theorem]{Corollary}
\theoremstyle{remark}
 \newtheorem{definition}[theorem]{Definition}
 \newtheorem{remark}[theorem]{Remark}
 \newtheorem*{remark*}{Remark}
 \newtheorem*{acknowledgements}{Acknowledgements}
 \newtheorem{example}[theorem]{Example}
\numberwithin{equation}{section}

\newcommand{\Z}{\boldsymbol{Z}}
\newcommand{\R}{\boldsymbol{R}}

\renewcommand{\phi}{\varphi}
\renewcommand{\epsilon}{\varepsilon}
\newcommand{\op}{\operatorname}

\newcommand{\pmt}[1]{{\begin{pmatrix} #1  \end{pmatrix}}}

\author{Shuntaro~Ohno}
\address[Ohno]{%
   Department of Mathematics, 
Kyoto Koka Senior high school,
Nishi-kyogoku Nodacho,
Kyoto 615-0861, Japan.
}
\email{tkmskr0329@yahoo.co.jp}

\author{Tetsuya~Ozawa}
\address[Ozawa]{%
Department of Mathematics, Meijo University,
Tempaku, Nagoya, 468-8502 Japan}
\email{ozawa@meijo-u.ac.jp}

\author{Masaaki Umehara}
\address[Umehara]{%
   Department of Mathematics, Graduate School of Science,
   Osaka University,
   Toyonaka, Osaka 560-0043,
   Japan
}
\email{umehara@math.sci.osaka-u.ac.jp}


\thanks{
The third author was partially supported by the Grant-in-Aid for 
Scientific Research (A) No.22244006, Japan Society for the 
Promotion of Science.}

\begin{document}

\begin{abstract}
We define a computable topological invariant $\mu(\gamma)$ for 
generic closed planar regular 
curves $\gamma$, which gives an effective lower bound
for the number of inflection points 
on a given generic closed planar curve.
Using it, we classify the topological types of 
\lq{\it locally convex curves}\rq\
(i.e. closed planar regular curves 
without inflections) whose numbers of crossings are 
less than or equal to five. 
Moreover, we discuss the relationship between
the number of double tangents and the invariant $\mu(\gamma)$
on a given $\gamma$.
\end{abstract}
\maketitle

%%%%%%%%%%%%%%%%%%%%%%%%%%%%%%%%%%%%%%%%%%%%%%%%%%%%%
\section{Introduction.}
\setcounter{section}{0}
In this paper, curves are always
assumed to be regular  (i.e. immersed).
The well-known Fabricius-Bjerre \cite{FB} 
theorem asserts (see also \cite{H1})
that
\begin{equation}\label{eq:FB}
d_1(\gamma)-d_2(\gamma)=\#_\gamma+\frac{i_\gamma}2
\end{equation}
holds for closed curves $\gamma$ satisfying 
suitable genericity assumptions,
where $d_1(\gamma)$ (resp. $d_2(\gamma)$) is 
the number of double tangents
of same side (resp. opposite side) and 
$\#_\gamma$ and $i_\gamma$ are the number of crossings 
and the number of inflections on $\gamma$, respectively.

\begin{figure}[htb]
\begin{center}
        \includegraphics[height=5.2cm]{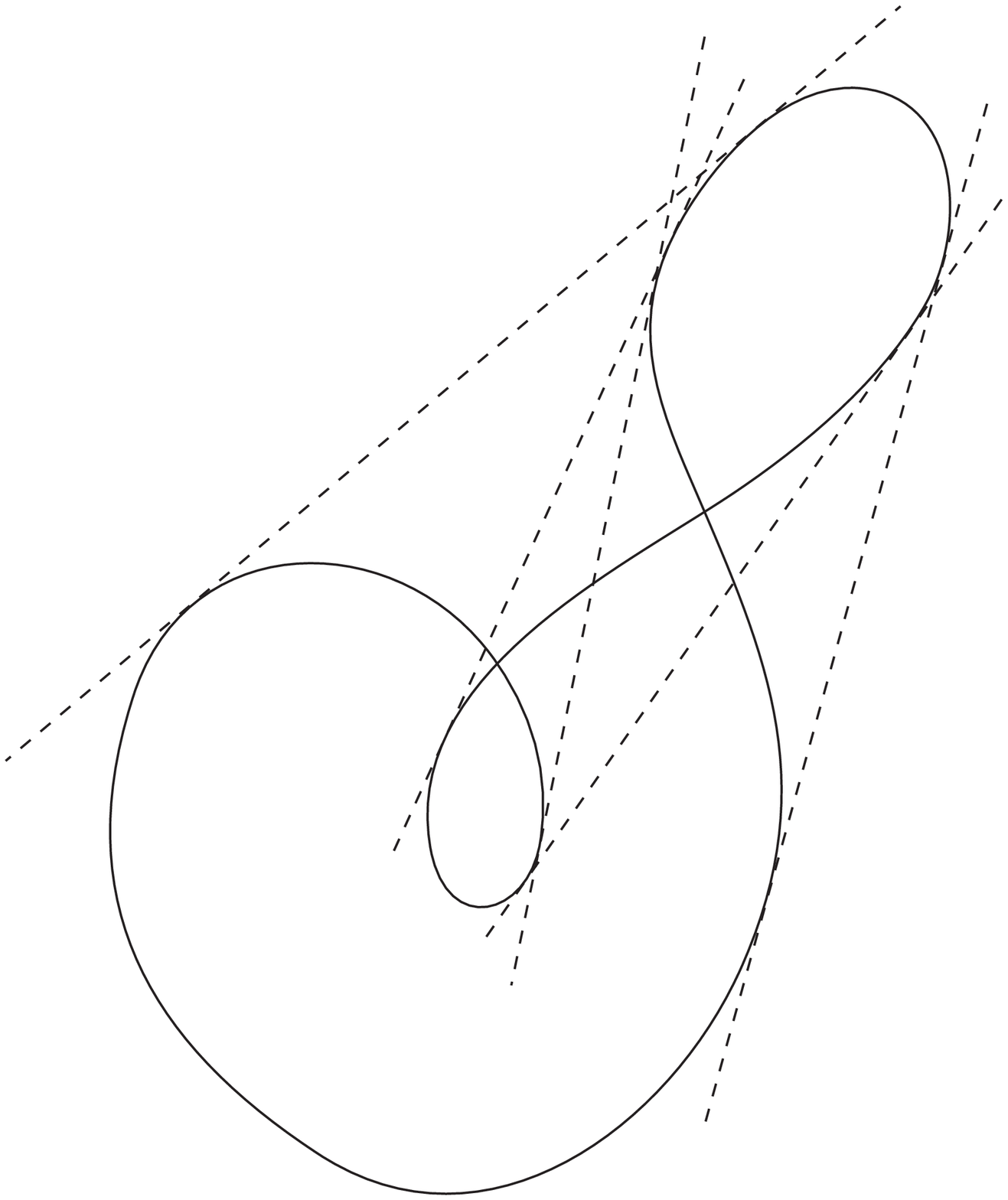}
\caption{A curve with $\#_\gamma=2$, $i_\gamma=2$, $d_1=4$ and $d_2=1$.}
\label{fig:fb}
\end{center}
\vspace{-2ex}
\end{figure}

However, for arbitrarily given  
non-negative integers $d_1,d_2,n$ and $i$
satisfying $d_1-d_2=n+i/2$, there might not exist
a corresponding curve, in general. 
As an affirmative answer to the Halpern conjecture in \cite{H2},
the second author \cite{O} proved the inequality 
\begin{equation}\label{eq:Ozawa}
d_1(\gamma)+d_2(\gamma)\le \#_\gamma(2\#_\gamma-1)
\end{equation}
for closed curves without inflections 
(see Remark \ref{rmk:halpern}).
It is then natural to expect that there might 
be further obstructions for
the topology of closed planar curves without 
inflections.
In this paper, we define a computable topological  invariant 
$\mu(\gamma)$ for closed planar curves.
By applying the Gauss-Bonnet formula, we show the 
inequality $i_\gamma\ge \mu(\gamma)$,
which is sharp at least for closed curves
satisfying $\#_\gamma\le 4$. In fact, for such $\gamma$,
there exists a closed curve $\sigma$ which  
has the same topological type
as $\gamma$ such that $i_\sigma= \mu(\sigma)$.
As an application, we classify the topological 
types of closed planar curves 
satisfying $i_\gamma=0$ and $\#_\gamma\le 5$. 
Moreover, we discuss the relationship between
the number of double tangents of the curve and the invariant $I(\gamma)$
on a given $\gamma$.

\section{Preliminaries and main results}
We denote by $\R^2$ the affine plane, and by $S^2$
the unit sphere in $\R^3$.
A closed curve $\gamma$ in $\R^2$ or $S^2$
is called {\it generic} if
\begin{enumerate}
\item all crossings are transversal, and
\item the zeroes of curvature are nondegenerate.
\end{enumerate}
By stereographic projection, 
we can recognize $S^2=\R^2\cup \{\infty\}$.
Two generic closed curves $\gamma_1$ and $\gamma_2$ 
in $\R^2$ (in $S^2$) are called {\it geotopic}, 
or said to have the {\it same topological type}, in $\R^2$ 
(resp. in $S^2$)
 if there is an orientation 
preserving diffeomorphism $\phi$ of $\R^2$ (resp. on $S^2$)
such that 
$
\op{Im}\gamma_2=\phi(\op{Im}\gamma_1). 
$
This induces an equivalence relation
on the set of closed curves.
We denote equivalence classes by 
$\langle \gamma_1 \rangle$ (resp. $[\gamma_1]$).
We fix a closed regular curve $\gamma:S^1\to \R^2$.
A point $c\in S^1$ is called an {\it inflection point}
of $\gamma$ if $\op{det}(\dot\gamma(t),\ddot\gamma(t))$ vanishes
at $t=c$. 
We denote by $i_\gamma$ the number of
inflection points on $\gamma$.
A closed curve  $\gamma:[0,1]\to \R^2$ 
is called {\it locally convex} 
if $i_\gamma=0$.
Whitney \cite{W} proved that 
any two closed  curves 
are regularly homotopic if and only
if their rotation indices coincide.
So then one can ask if this regular homotopy
preserves the locally convexity when the given
two curves are both locally convex.
In fact, one can easily prove this
by a modification of Whitney's argument,
which has been pointed out in \cite[p35, Exercise 11]{O2}.
For the sake of the readers' convenience we shall 
outline the proof:

\begin{proposition}
Let $\gamma_1$ and $\gamma_2$ be two 
locally convex closed regular curves.
Suppose that $\gamma_1$ has the same rotation 
index as $\gamma_2$.
Then there exists a family of closed 
curves 
$\{\Gamma_\epsilon\}_{\epsilon\in [0,1]}$
such that 
\begin{enumerate}
\item $\Gamma_0=\gamma_1$ and $\Gamma_1=\gamma_2$,
\item each $\Gamma_\epsilon$ $(\epsilon\in [0,1])$ is a locally convex
closed regular curve.
\end{enumerate}
\end{proposition}

\begin{proof}
Suppose that the curves $\gamma_j$ $(j=1,2)$
are both positively 
curved and have the same rotation indices, equal to $m$.
We change, if necessary, the parametrizations of $\gamma_j$
so that the tangent vectors $\dot\gamma_j(t)$ are 
positive scalar multiples of $\pmt{\cos t\\ \sin t}$.
Define a homotopy $\Gamma_\epsilon$ between $\gamma_1$ and 
$\gamma_2$ by
\[
  \Gamma_\epsilon(t) 
  = (1-\epsilon)\gamma_1(t) + \epsilon\gamma_2(t).
\]
It is easy to verify that 
$\{\Gamma_\epsilon\}_{\epsilon\in[0,1]}$ satisfies 
the required conditions (1) and (2).
\end{proof}

For a given generic closed curve $\gamma$,
we set
\begin{equation}
I(\gamma):=\min_{\sigma\in \langle \gamma\rangle}i_\sigma.
\end{equation}
Since $i_\gamma$ is an even number, so is $I(\gamma)$.
Inflection points on curves
in $\R^2$  correspond to
singular points on their Gauss maps. 
So it is natural to ask about the existence of
topological restrictions on closed curves 
without inflections, in other words,
we are interested in the topological
type of closed curves $\gamma$ satisfying 
$I(\gamma)=0$.

There are explicit combinatorial
procedures for determining generic
closed spherical curves with a 
given number of crossings,
as in Carter \cite{S} and Cairns and Elton \cite{CE}
(see also Arnold \cite{A}). 
We denote by $\#_\gamma$ the 
number of crossings for a given generic 
closed  curve. In this paper, we use the table 
of closed spherical curves
with $\#_\gamma\le 5$
given in the appendix of \cite{KU}.

\begin{figure}[htb]
\begin{center}
        \includegraphics[height=1.3cm]{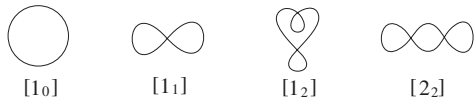}
\caption{Spherical closed curves with $\#_\gamma\le 2$.}
\label{fig:curve2s}
\end{center}
\vspace{-2ex}
\end{figure}

For example, the table of closed spherical curves 
with $\#_\gamma\le 2$ is given in Figure \ref{fig:curve2s},
where $1_2$ (resp. $2_2$)
means that the corresponding curve has 2-crossings and
appears in the table of curves in \cite{KU} with $\#_\gamma=2$ primary
(resp. secondary).

\begin{figure}[htb]
\begin{center}
        \includegraphics[height=2.8cm]{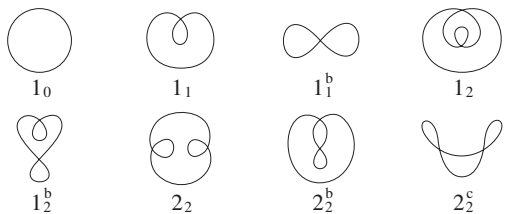}
\caption{Planar closed curves with $\#_{\gamma}\le 2$.}
\label{fig:curve2a}
\end{center}
\vspace{-2ex}
\end{figure}

Moving the position of $\infty$ 
via motions in $S^2=\R^2\cup \{\infty\}$,
we get the table of
closed planar curves with $\#_\gamma \le 2$
as in Figure \ref{fig:curve2a}. 
For example $1_1$ and $1_1^b$ 
(resp. $2_2$, $2_2^b$ and $2_2^c$) are
equivalent to $[1_1]$ (resp. $[2_2]$)
as spherical curves. Here, only the curves of
type $1_1$ and $2_2$
can be drawn with  no inflections.
Similarly, using the table of spherical
curves with $\#_\gamma\le 5$ given in \cite{KU},
we prove the following theorem.
The authors do not know of any reference for 
such a classification of generic locally convex curves.

\begin{theorem}\label{thm:main}
For a given generic closed regular curve $\gamma$
in $\R^2$, the inequality
\begin{equation}\label{eq:affine}
I(\gamma)\ge \mu(\gamma)
\end{equation}
holds. Moreover,
$I(\gamma)=0$  if and only if $\mu(\gamma)=0$ under the assumption
that $\#_\gamma\le 5$.
\end{theorem}

\begin{figure}[htb]
\begin{center}
        \includegraphics[height=1.2cm]{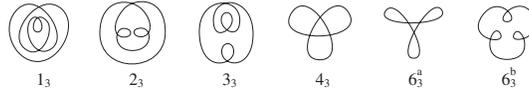}\\
\caption{Curves with $I_\gamma=0$ and $\#_\gamma= 3$ 
(6 curves in total).}
\label{fig:3}
\end{center}
\vspace{-2ex}
\end{figure}

\medskip
In particular, the number of equivalence classes of 
closed locally convex curves with $\#_\gamma\le 5$
is $76$ (see Figure \ref{fig:curve2a} and Figures
\ref{fig:3}, \ref{fig:4} and \ref{fig:5} in Section 3).
For example, in Figure \ref{fig:curve2a} 
the curves of type $1_0,1_1,1_2,2_2$ satisfy $I(\gamma)=0$,
and the remaining $1_1^b,1_2^b,2_2^b,2_2^c$ satisfy 
$I(\gamma)=2$. In Figure \ref{fig:3}, the curve of type $6_3^a$
is of the same topological type as $6_3^b$ as a spherical
curve, which is obtained from the 6th curve in
the table of curves with
 $\#_\gamma=2$ in the appendix of \cite{KU}.

\begin{figure}[htb]
\begin{center}
        \includegraphics[height=1.2cm]{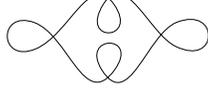}
\caption{A curve with $\#_\gamma=4$ and 
$I(\gamma) =\mu(\gamma) =4$.}
\label{fig:i4}
\end{center}
\vspace{-2ex}
\end{figure}

\begin{corollary}\label{cor:main}
A generic closed curve $\gamma$
with $\#_\gamma\le 4$ satisfies
$I(\gamma)= \mu(\gamma)\le 2$,
unless the topological type of 
$\gamma$ is as in Figure \ref{fig:i4}.
\end{corollary}

\begin{figure}[htb]
\begin{center}
        \includegraphics[height=2.5cm]{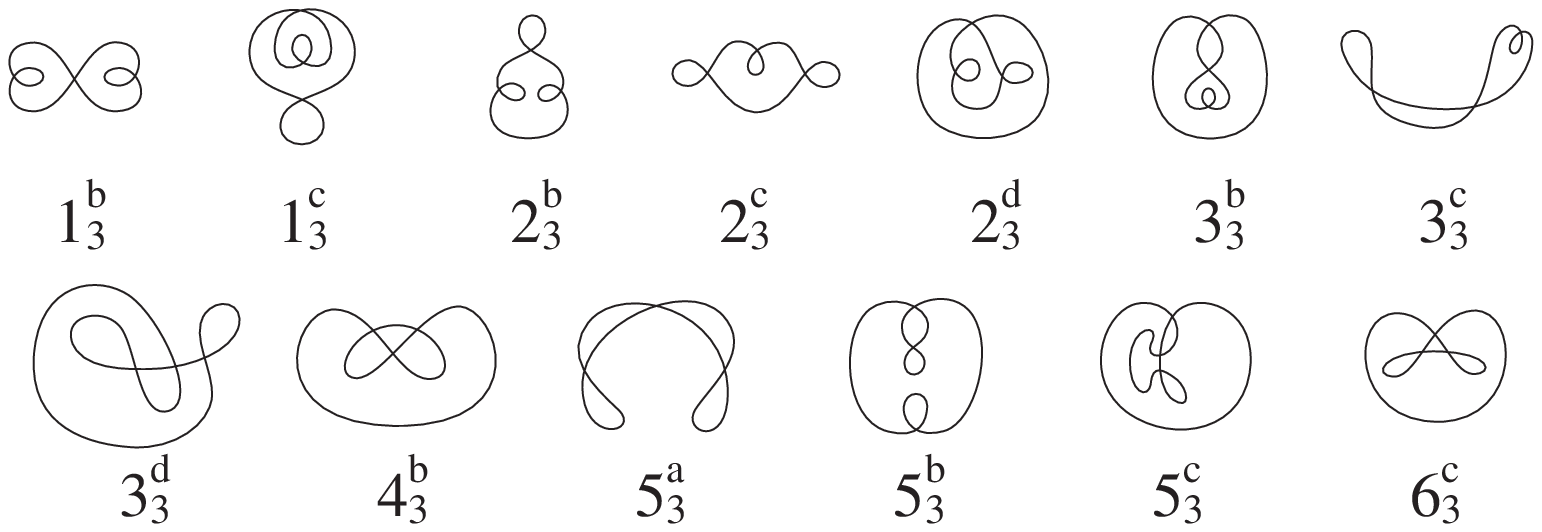}\qquad \quad 
\caption{Curves with $I(\gamma)>0$ and $\#_\gamma\le 3$ (13 curves 
in total).}
\label{fig:3g}
\end{center}
\vspace{-2ex}
\end{figure}

\noindent
A table of closed curves $\gamma$ with $i_\gamma=\mu(\gamma)$ 
for $\#_\gamma \le 3$ is given in 
Figures \ref{fig:curve2a}, \ref{fig:3} and \ref{fig:3g}.
For example, in Figure \ref{fig:3g}, 
the curves of type $6_3^c$ 
consideres as spherical curves
are of the same topological type as the curves 
of  type $6_3^a$ or $6_3^b$ 
in Figure \ref{fig:3}.

\section{Definition of the invariant $\mu(\gamma)$}

We fix a generic closed curve $\gamma:S^1\to \R^2$.
We set $\#_\gamma=m$.
We may suppose that $\gamma(0)=\gamma(1)$
is one of the crossings of $\gamma$.
Let 
$$
0=c_1<\cdots < c_{2m}(<1)
$$
be the inverse image of the crossings of
$\gamma$, which consists of $2m$ points in $S^1=\R/\Z$.
We set
$$
S^1_{\gamma}:=S^1\setminus \{c_1,...,c_{2m}\}.
$$
To introduce the invariant $\mu(\gamma)$,
we define special subsets on the curves called
{\lq $n$-gons\rq}:

\begin{figure}[htb]
\begin{center}
        \includegraphics[height=1.3cm]{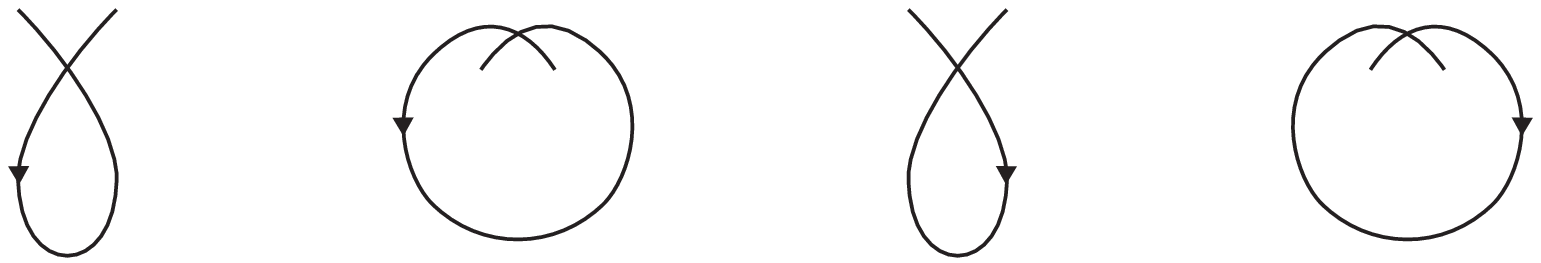}
\caption{Shells.}
\label{fig:shells}
\end{center}
\vspace{-2ex}
\end{figure}

\begin{definition}
Let $n(\ge 3)$ be an integer.
A disjoint union of $n$
proper closed intervals 
$$
J:=[a_1,b_1]\cup \cdots \cup [a_n,b_n]
$$
on $S^1$ is called an {\it $n$-gon}
if $a_1,b_1,\cdots,a_n,b_n\in \{c_1,...c_{2m}\}$
and the image $\gamma(J)$ is 
a piecewise smooth simple closed curve in $\R^2$.
The simply connected domain
bounded by $\gamma(J)$
is called the {\it interior domain} of the $n$-gon.
An $n$-gon is called {\it admissible} if 
at most two of the $n$ interior 
angles of $D$ are less than $\pi$.
\end{definition}

\begin{figure}[htb]
\begin{center}
        \includegraphics[height=1.3cm]{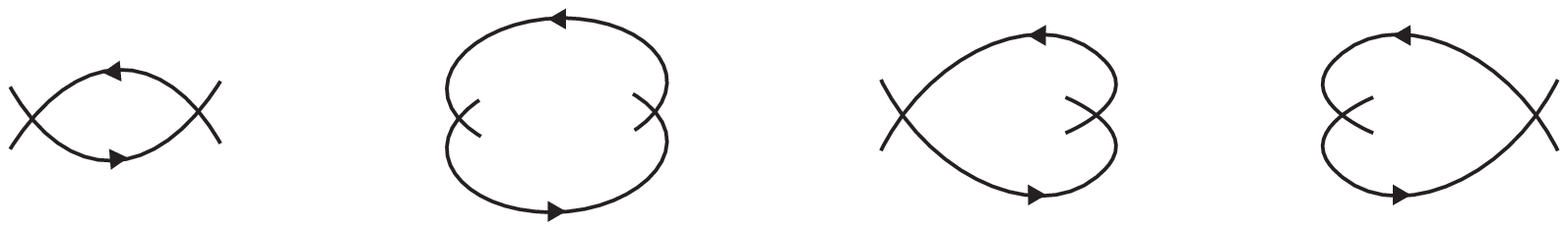}
\caption{Leaves.}
\label{fig:leaves}
\end{center}
\vspace{-2ex}
\end{figure}

\begin{figure}[htb]
\begin{center}
        \includegraphics[height=2cm]{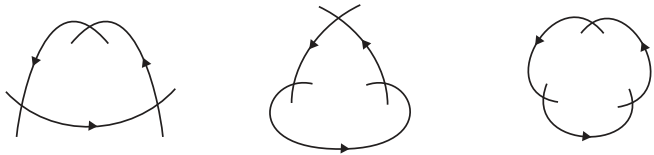}
\caption{Admissible triangles.}
\label{fig:triangles}
\end{center}
\vspace{-2ex}
\end{figure}

We denote by $\mathcal G_n(\gamma)$ the set of all 
admissible $n$-gons, and set
$$
\mathcal G(\gamma):=\bigcup_{n=1}^\infty \mathcal G_n(\gamma).
$$
Each element of $\mathcal G(\gamma)$ is called 
an {\it admissible polygon}.
A $1$-gon is called a {\it shell} 
(cf. Figure \ref{fig:shells}).
A $2$-gon is called a {\it leaf} 
(cf. Figure \ref{fig:leaves})
and a $3$-gon is called
a {\it triangle} (cf. Figure \ref{fig:triangles}).
All shells and all leaves are admissible.
However, a triangle whose interior angles are all
acute is not admissible.

We fix an admissible $n$-gon 
$
J:=[a_1,b_1]\cup \cdots \cup [a_n,b_n].
$
Then $\gamma(J)$ is a piecewise smooth simple closed curve
in $\R^2$. We give an orientation of $\gamma(J)$ so that
the interior domain of $\gamma(J)$ is on the 
left hand side of $\gamma(J)$.
This orientation induces an orientation on $[a_i,b_i]$ for each
$i=1,...,n$.
We call $[a_i,b_i]$  a {\it positive interval}
(resp. {\it negative interval}) 
if the orientation of the interval $[a_i,b_i]$ 
coincides with (resp. does not coincide with) the orientation of 
$\gamma$. 

Let $t$ be a point on $J\setminus\{a_1,b_1,...,a_n,b_n\}
=\bigcup_{i=1}^n (a_i,b_i)$.
Then $t$ belongs to an open interval $(a_i,b_i)$ for some
$i=1,...,n$. Then we set
$$
\epsilon_J(t):=
\begin{cases}
1 & \mbox{if $[a_i,b_i]$ is positive}, \\
-1 & \mbox{if $[a_i,b_i]$ is negative}. 
\end{cases}
$$

\begin{definition}
An admissible $n$-gon is called  {\it positive} 
(resp.  {\it negative})
if $\epsilon_J(t)>0$ (resp. $\epsilon_J(t)<0$)
holds for each 
$t\in J\setminus\{a_1,b_1,...,a_n,b_n\}$.
\end{definition}

The notions of positivity and negativity of shells were used
in \cite{U} and \cite{KU} differently from here.
The following assertion is the key to proving the
inequality \eqref{eq:affine}:

\begin{lemma} \label{lem:key}\ 
Let 
$
J:=[a_1,b_1]\cup \cdots \cup [a_n,b_n]
$
be an admissible $n$-gon.
Then there exists an interval $[a_i,b_i]$
{\rm ($1\le i\le n$)}
and  a point $c\in (a_i,b_i)$ such that
$\op{sgn}(\kappa_\gamma(c))=\epsilon_J(c)$, where
$
\kappa_\gamma(t):={\op{det}(\dot\gamma(t),\ddot \gamma(t))}/
{|\dot\gamma|^3}
$
is the curvature of $\gamma$
and $\op{sgn}(\kappa_\gamma(c))$ is the sign of the
real number $\kappa_\gamma(c)$.
\end{lemma}

\begin{proof}
Let $A_1,A_2,\cdots, A_n$ be the interior angles of 
the interior domain
of $J$, and set
$\Gamma=\gamma(J)$, which we regard as
an oriented piecewise smooth simple closed curve
with counterclockwise orientation.
In this proof,
$\R^2$ is considered as the Euclidean plane,
and we take the arclength parameter $s$ of $\Gamma$. 
Let $s=s_1,...,s_{n}$ be the points where $d\Gamma/ds$
is discontinuous.
We denote by $\kappa_{\Gamma}(s)$ 
($s\ne s_1,...,s_{n}$)
the curvature of the curve $\Gamma$.
Then, the Gauss-Bonnet formula yields that
$$
\int_{J}\kappa_{\Gamma}(s) ds=-(n-2)\pi+\sum_{i=1}^nA_i.
$$
Since $J$ is admissible, we may assume that $A_1,...,A_{n-2}\ge \pi$,
and so $\int_{J}\kappa_\Gamma(s) ds>0$.
Then there exist an index $i$ $(1\le i \le n)$ 
and $c\in [a_i,b_i]$
such that
$
\kappa_{\Gamma}(c)>0.
$
We denote by $\kappa_{\gamma}(s)$ the Gaussian curvature of 
$\gamma$ at $\Gamma(s)$.
Then we have that
$$
0<\kappa_{\Gamma}(c)=\epsilon_J(c)\kappa_{\gamma}(c),
$$
which proves the assertion.
\end{proof}

We now define the invariant $\mu(\gamma)$
mentioned in the introduction:

\begin{definition}
A function
$
\Phi:S^1_{\gamma}\to \{0,-1,1\}
$
is called an {\it admissible  
function} of $\gamma$
if it satisfies the following conditions:
\begin{enumerate}
\item 
$
\op{supp}(\Phi):=\Phi^{-1}(1)\cup \Phi^{-1}(-1)
$
is a finite set, and
\item for each $J\in \mathcal G(\gamma)$, 
there exists $t\in \op{supp}(\Phi)$
such that $t\in J$ and $\epsilon_J(t)=\Phi(t)$.
\end{enumerate}
A point $t\in \op{supp}(\Phi)$ is called a {\it positive point}
(resp. {\it negative point})
if $\Phi(t)>0$ (resp. $\Phi(t)<0$).
We denote by $\mathcal A$ the set of all admissible functions
of $\gamma$.
\end{definition}

We fix an admissible function $\Phi\in \mathcal A$.
Then we have an expression
$
\op{supp}(\Phi)=\{u_1,...,u_{\ell}\}
$
such that
$
c_1\le u_1< u_2<\cdots < u_{\ell}\le c_{2m}.
$
Let $\mu(\Phi)$ be the
number of sign changes
of the sequence
$$
\Phi(u_1),\,\,\Phi(u_2),\,\,\cdots,\Phi(u_{\ell}),\,\,\Phi(u_1).
$$
Then we set
\begin{equation}
\mu(\gamma):=\min_{\Phi\in \mathcal A}\mu(\Phi).
\end{equation}

\medskip
\noindent
{\it Proof of the inequality \eqref{eq:affine}.}
Let $\gamma$ be a generic closed curve in $\R^2$.
We take a curve $\sigma\in \langle \gamma\rangle$
such that $i_\sigma=I(\gamma)$.
Without loss of generality, we may assume that
$i_\gamma=I(\gamma)$.
We can take a point $c_J\in [0,1]$ for each 
$J\in \mathcal G$ in order that 
$\op{sgn}(\kappa_\gamma(c_J))=\epsilon_J(c_J)$,
and that $c_{J}\neq c_K$ if $J\neq K$.
We define a function $\Phi:S^1_\gamma\to \{0,-1,1\}$
by
$$
\Phi(t):=
\begin{cases}
\op{sign}(\kappa_\gamma(c_J)) & \mbox{if $t=c_J$ for some 
$J\in \mathcal G(\gamma)$},
 \\
0 & \mbox{otherwise}.
\end{cases}
$$
Then, $\Phi$ is an admissible  
function.
Since $i_\gamma=I(\gamma)$,
the curvature function of $\gamma$ changes sign at most 
$I(\gamma)$ times. So we have that
$\mu(\Phi)\le I_\gamma$,
in particular, $\mu(\gamma)\le I_\gamma$.
\qed

Also we have the following assertion:

\begin{proposition}
\label{prop:trivial}
Let $\gamma$ be a generic closed curve in $\R^2$.
Then $\mu(\gamma)$ is  a non-negative even integer,
as well as $I(\gamma)$. Moreover, 
$\mu(\gamma)>0$ holds if and only if
$\mathcal G(\gamma)$ does not contain a positive
polygon and a negative polygon at the same time.
\end{proposition}

\begin{proof}
Since the number of sign changes of a cyclic sequence 
of real numbers is always even, $\mu(\gamma)$ is also even.
Moreover,
if $\gamma$ has a positive (resp. negative)
polygon, each admissible 
function
$\Phi$ must take a positive (resp. negative)
value. 
So the existence of two distinct polygons of
opposite sign implies that $\mu(\gamma)\ge 2$.

Now, we prove the converse.
A closed curve which is not a simple closed curve $1_0$
has at least one shell, and a shell is necessarily 
a positive or a negative polygon. Suppose that $\gamma$
has no negative polygons.
Then we can choose a point $c_J\in J\cap S^1_\gamma$ 
for each admissible polygon $J\in \mathcal G(\gamma)$
such that
$\epsilon_J(c_J)>0$.
If we set
$$
\Phi(t):=
\begin{cases}
\epsilon_J(c_J) & 
\mbox{if $t=c_J$ for some $J\in \mathcal G(\gamma)$}, \\
0 & \mbox{otherwise},
\end{cases}
$$
then $\Phi$ is an admissible function, and
$\mu(\Phi)=0$. Thus we have $\mu(\gamma)=0$.
This proves the converse.
\end{proof}

To show the computability of the invariant $\mu(\gamma)$,
we fix $4m$ ($m=\#_\gamma$)
points $t_i,s_i$ ($i=1,...,2m$) satisfying
$$
0=c_1<t_1<s_1<c_2<\cdots < c_{2m}<t_{2m}<s_{2m}(<1),
$$
and show the following lemma:

\begin{lemma}\label{lemma:reduction}
For each admissible function $\Phi$ of $\gamma$, 
there exists
an admissible function $\Psi$ satisfying the 
following properties:
\begin{enumerate}
\item $\op{supp}(\Psi)\subset \{t_1,s_1,...,t_{2m},s_{2m}\}$,
\item $\mu(\Psi)\le  \mu(\Phi)$.
\end{enumerate}
\end{lemma}

\begin{proof}
We fix an interval $U=(c_i,c_{i+1})$, where
$c_{2m+1}:=c_1$.
If $\Phi$ takes non-negative (resp. non-positive)
values on $U$,
then we set
$$
\Psi(t):=
\begin{cases}
1 \mbox{\,\,(resp. $-1$)} & \mbox{if $t=t_i$}, \\
0 & \mbox{if $t\in U\setminus\{t_i\}$}.
\end{cases}
$$
If $\op{supp}(\Phi)\cap U$ contains two points
$v_1,v_2$ such that $v_1<v_2$ and
$\Phi(v_1)=-\Phi(v_2)$, then we set
$$
\Psi(t):=
\begin{cases}
\Phi(v_1) & \mbox{if $t=t_i$}, \\
\Phi(v_2) & \mbox{if $t=s_i$}, \\
0 & \mbox{otherwise}.
\end{cases}
$$
Since $U$ is arbitrary, we get
a function $\Psi$ defined on $S^1_\gamma$.
Since $U \cap J$ coincides with
 either $U$ or an empty set
for each $J\in \mathcal G(\gamma)$,
$\Psi$ is an admissible function, and one can easily verify
$\mu(\Psi)\le \mu(\Phi)$.
\end{proof}

\begin{remark}\label{rmk:computability}
(Computability of the invariant $\mu(\gamma)$)
The function $\Psi$ obtained in Lemma \ref{lemma:reduction}
is called a {\it reduction} of $\Phi$.
($\Psi$ may not be uniquely determined from $\Phi$,
since 
$\op{supp}(\Phi)\cap U$ might consist of more than two points.)
By definition, there exists an admissible  
function $\Phi$
such that $\mu(\Phi)=\mu(\gamma)$.
By Lemma \ref{lemma:reduction}, there is a
reduction of such a function $\Phi$. 
Thus the invariant $\mu(\gamma)$ is attained by
a reduced function $\Psi$.
Since the number of reduced admissible  
functions
is at most $3^{4m}$, the invariant $\mu(\gamma)$ can be
computed in a finite number of steps.
\end{remark}

\begin{remark}\label{rmk:flexibility}
(A flexibility of the reduced admissible  
function)
In the above construction of the function $\Psi$ via
$\Phi$ we may set
$$
\Psi(t):=
\begin{cases}
\Phi(v_2) & \mbox{if $t=t_i$}, \\
\Phi(v_1) & \mbox{if $t=s_i$}, \\
0 & \mbox{otherwise},
\end{cases}
$$
when $\op{supp}(\Phi)\cap U$ contains two points
$v_1,v_2$ such that $v_1<v_2$ and
$\Phi(v_1)=-\Phi(v_2)$.
Then $\Psi$ is also an admissible  
function.
This modification of $\Psi$ can be done for each 
fixed interval $U=(c_i,c_{i+1})$.
However, after the operation, it might not hold that
$\mu(\Psi)\le  \mu(\Phi)$.
\end{remark}

\begin{proposition}
\label{prop:crossings}
Let $\gamma$ be a generic closed curve in $\R^2$.
Then it satisfies
\begin{equation}\label{eq:hidden}
\mu(\gamma)\le 2\#_\gamma.
\end{equation}
\end{proposition}

\begin{proof}
As pointed out in Remark \ref{rmk:computability},
there exists a reduced admissible function $\Psi$.
Since $\mu(\Psi)\le 4m$, we get the estimate 
$\mu(\gamma)\le 4\#_\gamma$. 
However, we can improve it as follows:
As pointed out in Remark \ref{rmk:flexibility},
one can replace the values $\Psi(t_i),\Psi(s_i)$
by $-\Psi(t_i),-\Psi(s_i)$
whenever $\Psi(t_i)=-\Psi(s_i)$ for each $i=1,...,2m$. 
So using this modification inductively for $i=1,...,2m$
if necessary,
we can modify $\Psi$ so that $\mu(\Psi)\le 2m$, and 
then we have $\mu(\gamma)\le 2\#(\gamma)$.
\end{proof}

\begin{remark}
If $\gamma$ is a lemniscate $1^b_1$, then 
$\mu(\gamma)= 2\#_\gamma=2$ holds.
However, the authors do not know of any other example
satisfying the equality $\mu(\gamma) = 2\#_\gamma$
(cf. Question 3).
\end{remark}

\begin{example}\label{ex:c2}(Curves with 
a small number of intersections.)
Here, we demonstrate how to determine $\mu(\gamma)$
with $\#_\gamma\le 2$. In the eight classes of curves
in Figure \ref{fig:curve2a},
four classes have been drawn without inflections.
So $I(\gamma)=0$ holds for these four curves.
Each of the remaining four curves satisfies $\mu(\gamma)>0$,
by Proposition \ref{prop:trivial}.
On the other hand, these four curves 
in Figure \ref{fig:curve2a}
have been drawn with exactly two inflections.
Thus we can conclude that they satisfy $\mu(\gamma)=2$.
Now, let $\gamma$ be a curve 
as in Figure \ref{fig:i4}.
Then $\gamma$ has four disjoint shells, 
two of which are positive, and the other two
are negative. So we can conclude that
$\mu(\gamma)\ge 4$. Since $\gamma$ as in
Figure \ref{fig:i4} has exactly four inflections,
we can conclude that $I(\gamma)=\mu(\gamma)=4$.
\end{example}

\begin{figure}[htb]
\begin{center}
        \includegraphics[width=3.3cm]{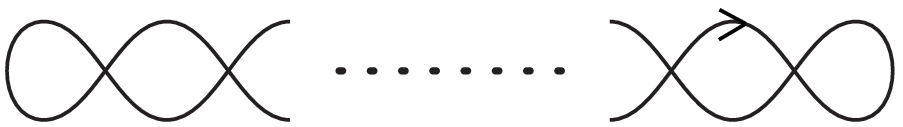}\qquad \quad 
        \includegraphics[width=2.8cm]{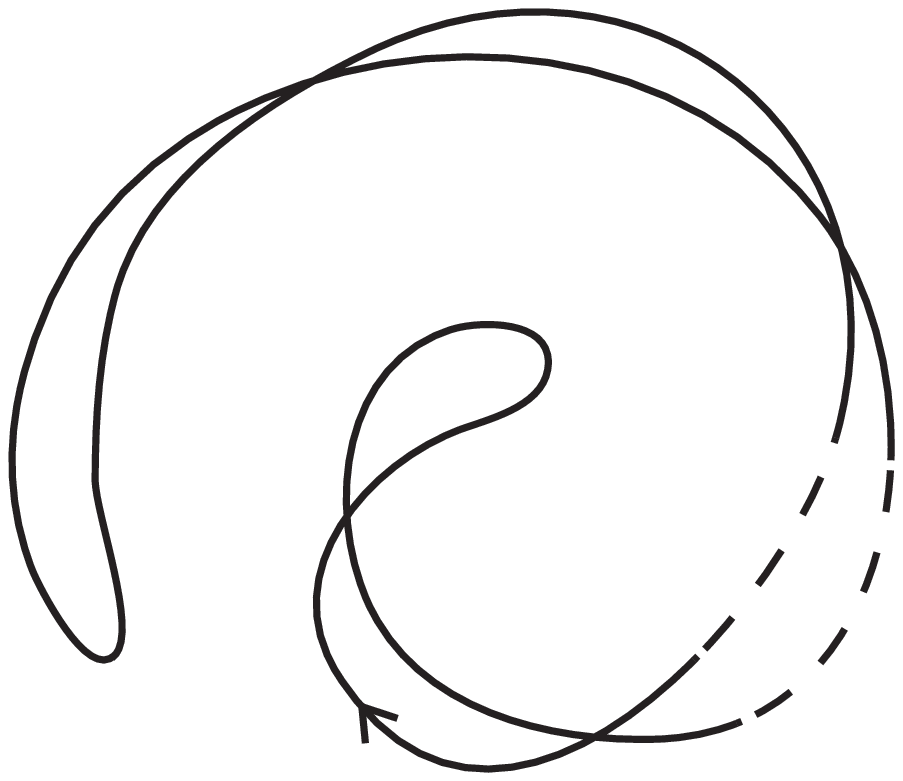}
\caption{A chain-like curve and its realization with 
$i_\gamma=2$.}
\label{fig:long}
\end{center}
\vspace{-2ex}
\end{figure}

\begin{example}\label{ex:chain} (Chain-like curves.)
We consider a curve $\gamma$ with $\#_\gamma=n$
($n\ge 1$) as in Figure \ref{fig:long}, left.
This curve has two shells and $n-1$ leaves,
including a positive shell and a negative leaf,
which are disjoint. 
Thus $(I(\gamma)\ge )\mu(\gamma) \ge 2$ holds.
As in Figure \ref{fig:long}, right,
this curve can be drawn along a spiral with 
two inflections. So we can conclude that
$I(\gamma)=\mu(\gamma)=2$.
In this manner, drawing curves along a spiral
is often useful to reduce the number of inflection points.
Several useful techniques for drawing curves 
with a restricted number of inflections
are mentioned in Halpern \cite[Section 4]{H2}.
\end{example}

\begin{figure}[htb]
\begin{center}
        \includegraphics[height=3.0cm]{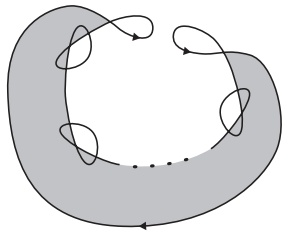} 

\caption{A figure of the curve  $\gamma_n$.}
\label{fig:n}
\end{center}
\vspace{-2ex}
\end{figure}

\begin{figure}[htb]
\begin{center}
        \includegraphics[height=2.3cm]{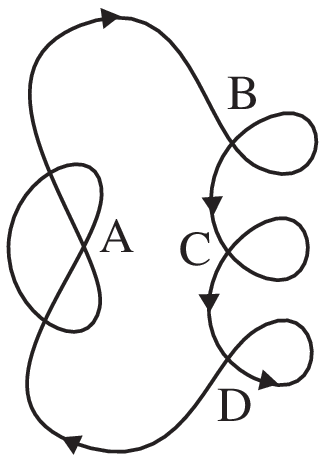}\qquad \qquad \quad 
        \includegraphics[height=2.5cm]{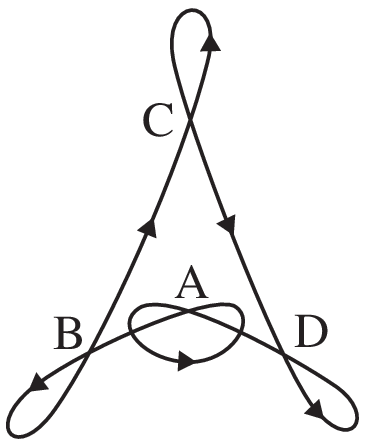}

\caption{A curve satisfying  $i_\gamma=0$.}
\label{fig:okey}
\end{center}
\vspace{-2ex}
\end{figure}

\begin{example}\label{ex:3} 
(Curves with a negative $n$-gon.)
We consider a curve $\gamma_n$ 
as in 
Figure \ref{fig:n},
which has several positive 
polygons, but only
one negative admissible  $n$-gon, 
marked in gray in Figure \ref{fig:n}.
So we can conclude that $I(\gamma)=\mu(\gamma)=2$,
as in Figure \ref{fig:n}.
This example shows that an $n$-gon 
$(n\ge 4)$ is needed
to find a curve that cannot be
locally convex.

Admissibility of a polygon is important for 
the definition of the invariant $\mu$.
In fact, the curve 
as in Figure \ref{fig:okey}, left, has 
negative polygons which are not admissible, and
it can be realized 
without inflections as in Figure \ref{fig:okey}, right.
\end{example}

\begin{figure}[htb]
\begin{center}
        \includegraphics[height=3.3cm]{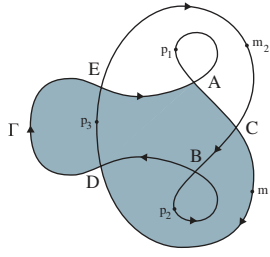}\qquad \qquad \quad 
\caption{A curve satisfying  $i_\gamma=4$ and $\#_\gamma=5$.}
\label{fig:ohno}
\end{center}
\vspace{-2ex}
\end{figure}

\begin{example}\label{ex:ohno} 
(A curve with an effective leaf which 
is neither positive nor negative.)
We consider a curve $\gamma$ as 
in Figure \ref{fig:ohno}. 
This curve has $5$ crossings
and exactly two positive shells at $A$ and $B$.
It also has a negative shell at $C$.
Thus $\mu\ge 2$ by Proposition \ref{prop:trivial}.
We show that $\mu(\gamma)=4$ by way of contradiction:
There exist two positive points $p_1$ and $p_2$ on 
the two positive shells at $A$ and $B$, respectively.
There is a unique simple closed arc $\Gamma$ bounded by 
$A$ and $B$ which passes through $D$ and $E$.
Suppose that $\mu(\gamma)=2$. Then there are no
negative points on $\Gamma$.
Now we look at the negative leaf with vertices  $A$ and $D$.
In Figure \ref{fig:ohno}, this leaf is marked in gray.
Since $\Gamma$ does not contain a negative point,
there must be a negative point $m_1$ between $A$ 
and $D$ on this leaf.
Since the curve has a symmetry,
applying the same argument to
the negative leaf at $B$ and $E$,
there is another negative point $m_2$
between $E$ and $C$.

Finally, we look at a leaf with vertices $D$ and $E$, which is
not positive nor negative.
Since $\Gamma$ has no negative point, there must have
a positive point $p_3$ on the arc on the 
right-hand side of the leaf.
Since the sequence
$p_1,m_1,p_3,m_2,p_2$ changes sign four times, 
this gives a contradiction. Thus 
$\mu(\gamma)\ge 4$.
Since the curve can be drawn with exactly four inflections
as in Figure \ref{fig:ohno},
we can conclude that $I(\gamma)=\mu(\gamma)=4$.
In this example, an admissible polygon which 
is neither positive nor negative plays
a crucial role to estimate the 
invariant $I(\gamma)$
by using $\mu(\gamma)$.
\end{example}

\begin{figure}[htb]
\begin{center}
  \includegraphics[height=3.5cm]{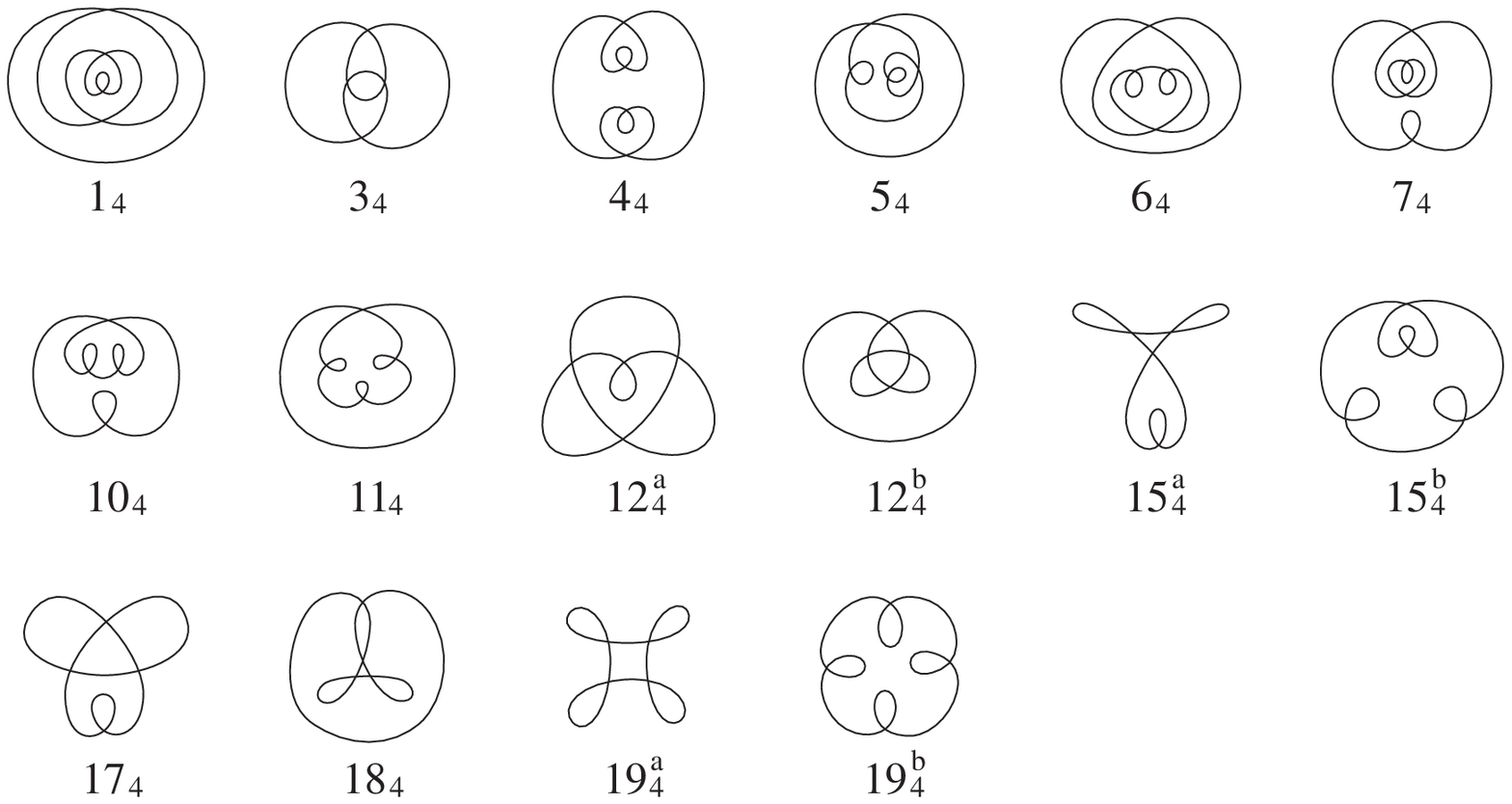} 
\caption{Curves with $I(\gamma)=0$ and $\#_\gamma\le 4$ (16 curves 
in total).}
\label{fig:4}
\end{center}
\vspace{-2ex}
\end{figure}

\begin{figure}[htb]
\begin{center}
        \includegraphics[height=6.5cm]{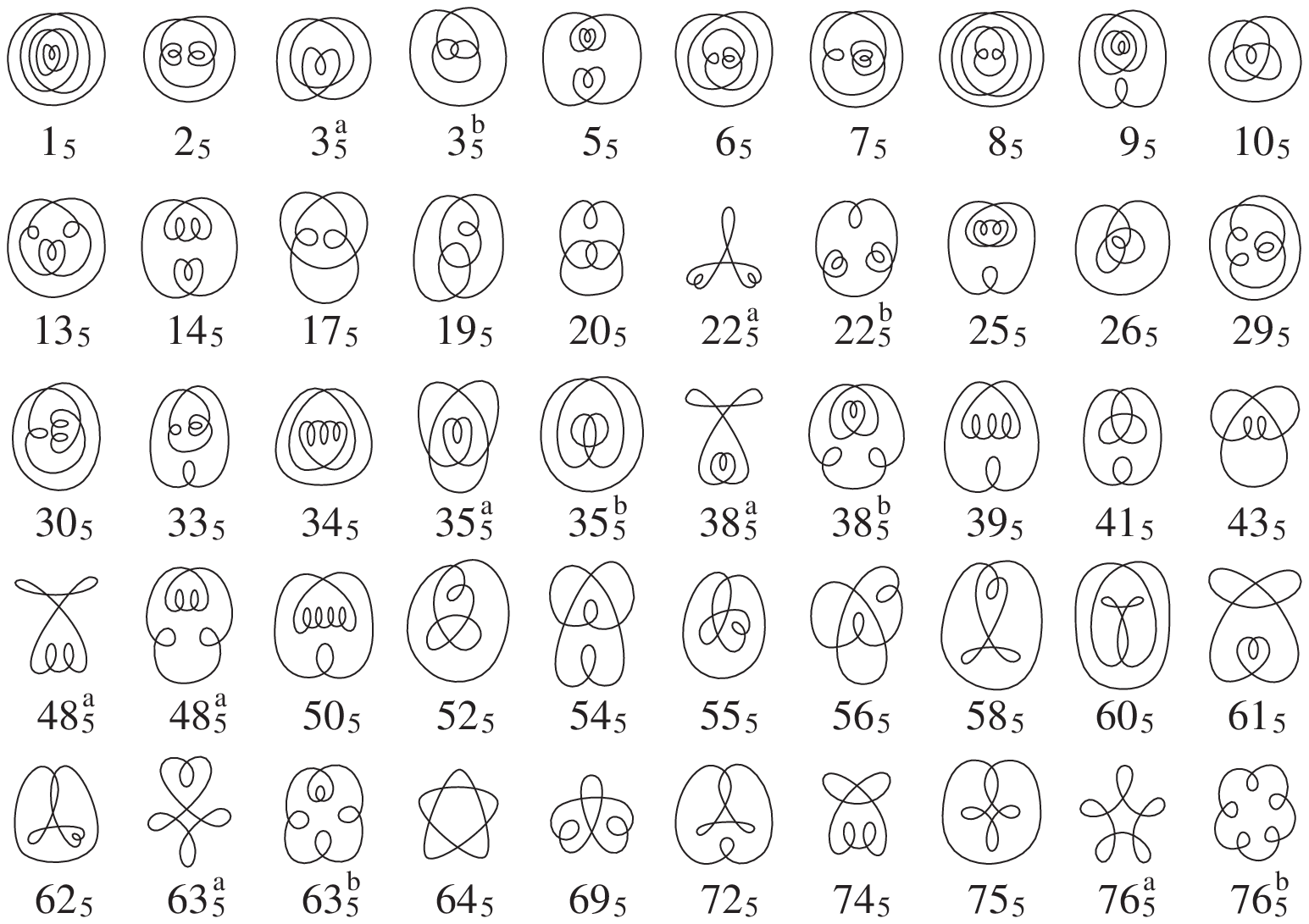}\qquad \quad 
\caption{Curves with $I(\gamma)=0$ and $\#_\gamma=5$ (50 curves
in total).}
\label{fig:5}
\end{center}
\vspace{-2ex}
\end{figure}

\medskip
\noindent
({\it Proof of the second assertion of Theorem \ref{thm:main}.})
The table of spherical curves up to $\#_\gamma\le 5$
is given in the appendix of \cite{KU}.
By moving the position of $\infty$ in $S^2=\R^2\cup\{\infty\}$,
we get the table of planar curves up to $\#_\gamma\le 5$
and can compute the invariant $\mu(\gamma)$.
So we can list the curves with $\mu(\gamma)=0$.
By Proposition \ref{prop:trivial}, it is sufficient 
to check for
the existence of positive polygons and negative polygons.
When $\#_\gamma\le 2$, the number of topological
types of such  curves is $3$.
If $\#_\gamma= 3,4,5$, then the number of 
topological types of such curves is $6,16,50$, respectively.
After that we can draw the pictures of 
the curves by hand. 
If we are able to draw $76$ figures of the curves without inflections,
the proof is finished, 
and this was accomplished in Figures \ref{fig:3},
\ref{fig:4} and \ref{fig:5}.
\qed

\medskip
\noindent
({\it Proof of Corollary \ref{cor:main}.})
For curves $\gamma$ with $\mu(\gamma)=0$, 
we can show $I(\gamma)=0$ by drawing curves
without inflections.
On the other hand, when $\mu(\gamma)>0$,
we can show $I(\gamma)=2$ by drawing curves
with exactly two inflections, except for 
the curve as in Figure \ref{fig:i4}.
\qed

\section{Double tangents and geotopical tightness}
In this section, we would like to give an application.

\begin{figure}[htb]
\begin{center}
        \includegraphics[height=3.5cm]{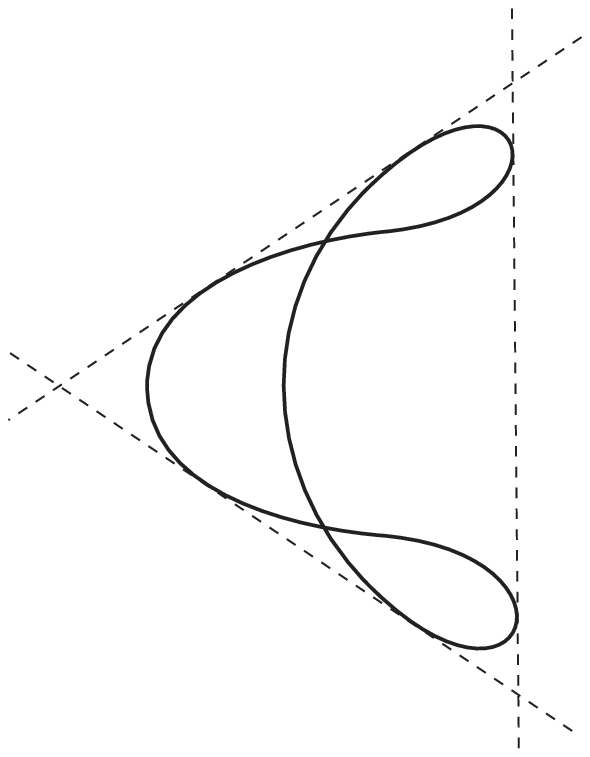}\qquad \qquad \quad 
        \includegraphics[height=4.0cm]{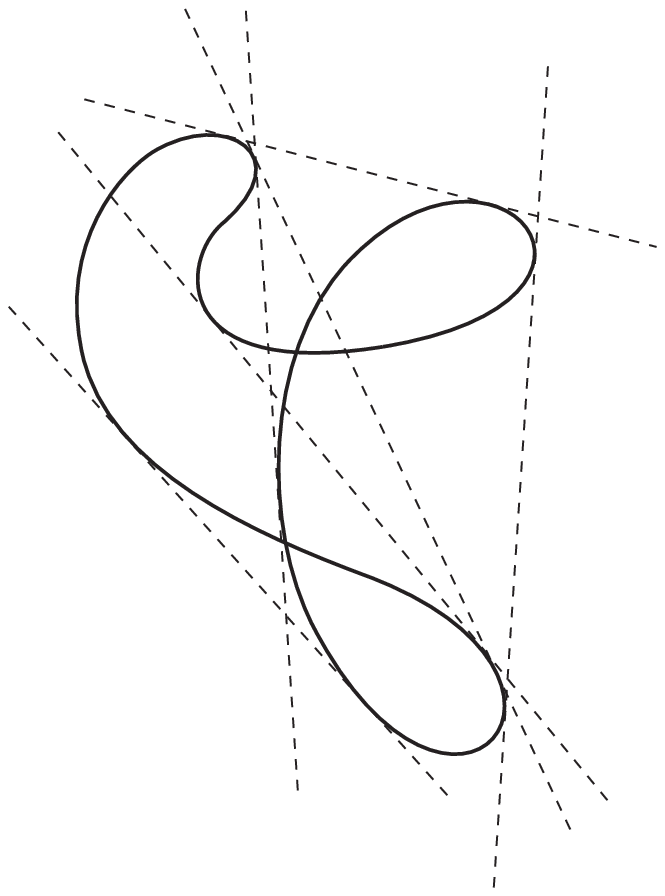}
\caption{Curves of $2_2^c$ with $i_\gamma=2$ 
but with different $d$.}
\label{fig:fb2}
\end{center}
\vspace{-2ex}
\end{figure}

\begin{definition}
Let $\gamma$ be a generic planar curve.
We set
\begin{equation}
d(\gamma):=d_1(\gamma)+d_2(\gamma).
\end{equation}
Then we define 
\begin{equation}
\delta(\gamma):=\min_{\sigma\in \langle \gamma\rangle}d(\sigma),
\end{equation}
which gives the minimum number of double tangents
of  the curves in the equivalence class $\langle \gamma\rangle$.
(As in Figure \ref{fig:fb2}, $d(\sigma)$ might be different from
$d(\gamma)$ even if $\sigma\in \langle \gamma\rangle$ and $i_\gamma
=i_\sigma=I(\gamma)$.) A curve $\sigma\in \langle \gamma\rangle$ satisfying 
$d(\gamma)=\delta(\gamma)$
is called {\it geotopically tight} or {\it $g$-tight}.
We call the integer $\delta(\gamma)$ 
the {\it $g$-tightness number}.
\end{definition}

The following assertion holds:

\begin{proposition}\label{prop:d}
It holds that
$\delta(\gamma)\ge \#_\gamma+I(\gamma)/2$.
\end{proposition}

\begin{proof}
By \eqref{eq:FB}, we get that
$$
d(\sigma)=d_1(\sigma)+d_2(\sigma)
\ge d_1(\sigma)-d_2(\sigma)=\#_\sigma+\frac{I(\sigma)}2
$$
holds for $\sigma\in \langle \gamma\rangle$.
Taking the infimum, we get the assertion.
\end{proof}

We expect that any locally convex
generic closed curves might be $g$-tight
(see Question 4).
Relating this, we can prove the following
assertion:

Let $\gamma:\R\to\R^2$ be a periodic parametrization 
by arclength of a locally convex curve with total 
length $\ell$, that is,
 $\gamma(t+\ell)=\gamma(t)$ holds for $t\in\R$. 
Without loss of generality, we may assume 
that the curvature of $\gamma$ is positive. 
Let 
$$
x=\gamma(s_1)=\gamma(s_2)\qquad (0<|s_1-s_2|<\ell)
$$
be a crossing of $\gamma$.
Replacing $s_1$ by $s_1+\ell$ if necessary, 
we may assume that $s_1<s_2$ and
$\dot\gamma(s_2)$ is a positively rotated vector of 
$\dot\gamma(s_1)$ through an angle $\alpha$ with $0<\alpha<\pi$. 
When the parameter $t$ varies in the interval $[s_1,s_2]$, 
the tangent vector $\dot\gamma(t)$ 
rotates through an angle $2\pi n_1+\alpha$,
where $n_1$ is a positive integer.
Similarly, for the interval $[s_2,s_1+\ell]$, 
there exists a positive integer $n_2$
such that the rotation angle
of $\dot\gamma(t)$ is equal to $2\pi n_2-\alpha$.
The sum $n_1+n_2$ is the total rotation index of $\gamma$.
Denote by $W(x)$ the difference $n_1-n_2$.
We easily recognize that the sum of $W(x)$ for each crossing
$x$ of $\gamma$ is a geotopy invariant.
The following theorem can be proved using the equality
of $d_2$ in \cite[p7]{O} (we omit the details):

\begin{theorem}\label{prop:conj}
The number of double tangents  for
any locally convex generic closed
curve $\gamma$ depends only on its 
geotopy type. More precisely, the following identity
holds:
\begin{equation}\label{eq:new}
d_2(\gamma)
=\sum W(x),
\end{equation}
where the sum runs over all crossings $x$ of $\gamma$.
\end{theorem}

\begin{remark}\label{rmk:halpern}
The rotation index $R_\gamma$ of $\gamma$ which is
at each crossing equal to $n_1+n_2$, as mentioned above,
is less than or equal to 
$\#_\gamma+1$ (cf. [{W}]).
Thus the formula \eqref{eq:new} implies
$$
d_2=\#_\gamma R_\gamma-2\sum W(x)\le
\#_\gamma (\#_\gamma+1)-2\#_\gamma=\#_\gamma (\#_\gamma-1),
$$
which reproves Halpern's conjecture in [{H2}].
\end{remark}

\begin{corollary}\label{cor:d}
$\delta(\gamma)=1,2,2,3,3$ for $\gamma$ of
type $1_1,1_1^b,1_2,2_2^b,2_2^c$ as in Figure \ref{fig:curve2a},
respectively.
\end{corollary}

\begin{proof}
By Proposition \ref{prop:d},
$\delta(\gamma)\ge 1,2,3,3$ for $\gamma$ of
type $1^b,2_2^b,2_2^c$, respectively.
On the other hand, the curves given 
in Figure \ref{fig:curve2a} attain equality
in this inequality.
\end{proof}

In Figure \ref{fig:curve2a}, there are two remaining
curves of type $1_2^b$ and $2_2$, 
whose g-tightness numbers have not been specified
by the authors.
The curve of type  $2_2$ given in Figure \ref{fig:curve2a}
satisfies $d(\gamma)=4$, and we expect that $\delta(\gamma)=4$.
On the other hand,
the corresponding curve as in 
the right of Figure \ref{fig:curve2a} 
satisfies $d=7$, but one can realize the curve
with $d=5$ as in Figure \ref{fig:fb}.
We expect that it might be
$g$-tight. If true, $\delta(\gamma)=5$ holds.

\begin{figure}[htb]
\begin{center}
        \includegraphics[height=3cm]{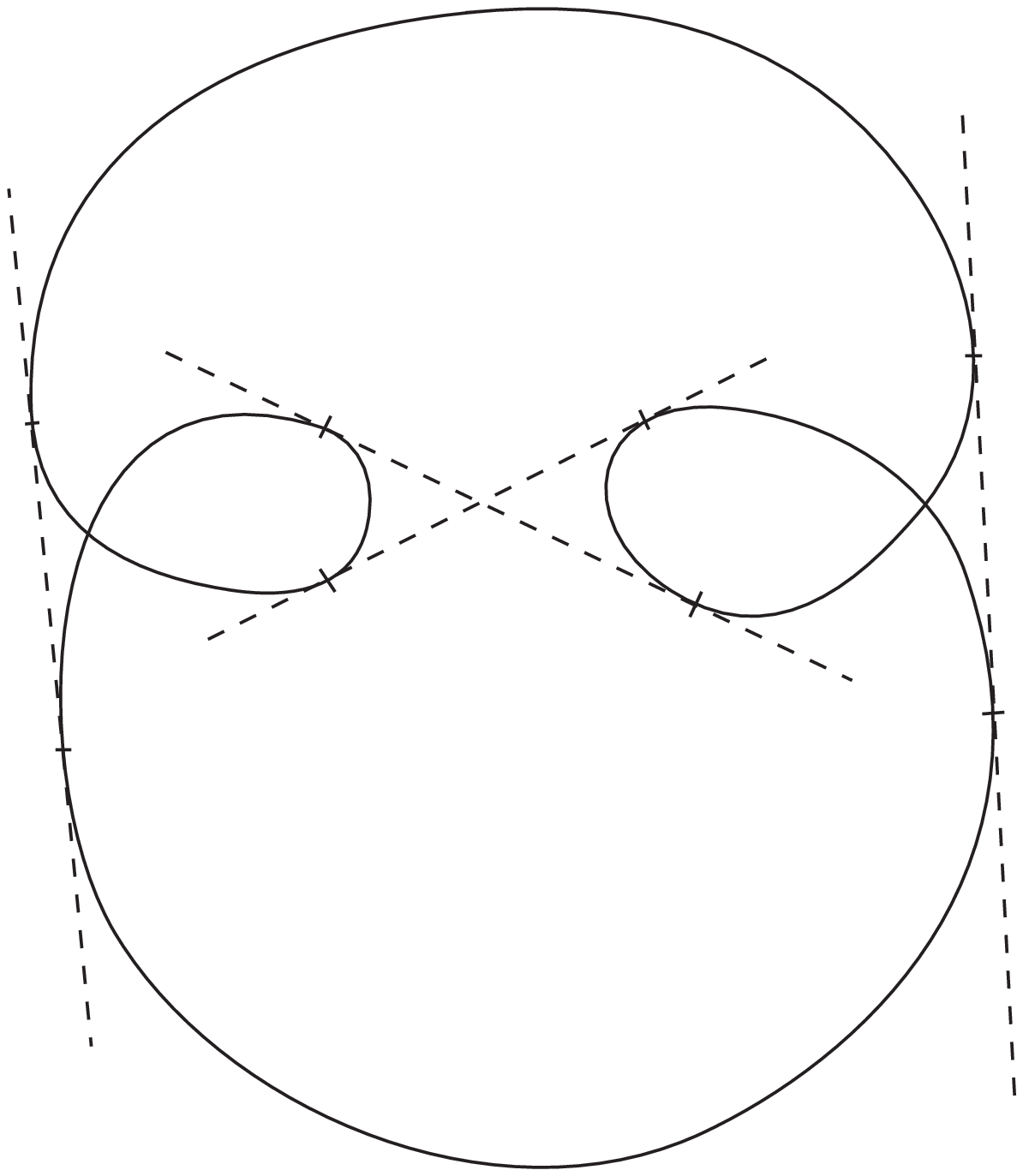}\qquad \qquad \quad 
        \includegraphics[height=3cm]{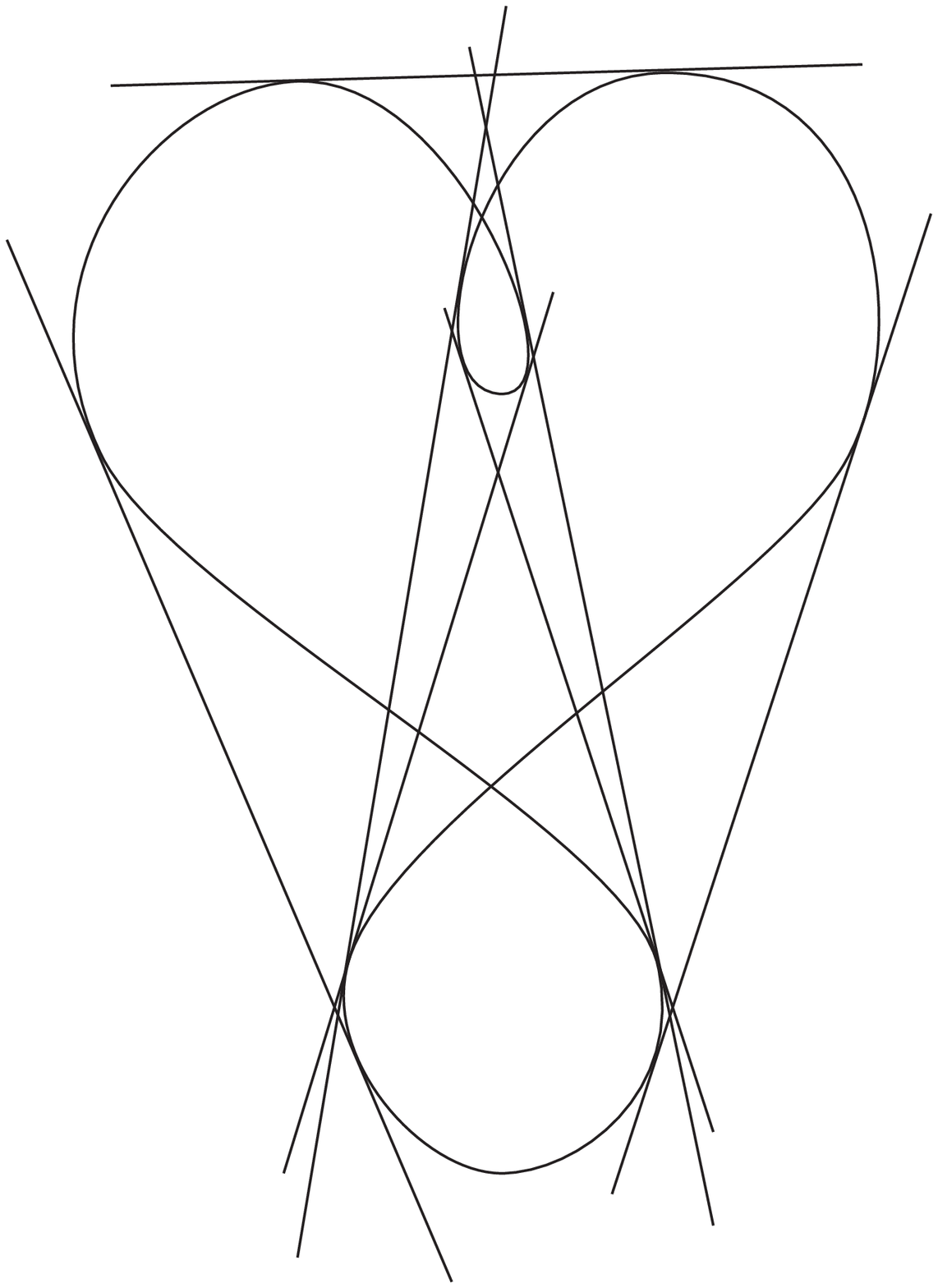}
\caption{Curves of type $2_2$ and $1^b_2$.}
\label{fig:75}
\end{center}
\vspace{-2ex}
\end{figure}

\bigskip
A relationship between inflection points and
double tangents for simple closed curves 
in the real projective plane
with a suitable convexity is given in \cite{TU}.
Finally, we leave several open questions on
the invariants $I(\gamma)$ and $\delta(\gamma)$:

\medskip
\noindent
{Question 1.}
{\it Does $\mu(\gamma)=0$ imply
$I(\gamma)=0$?}

\medskip
\noindent
{Question 2.}
{\it Is there a generic closed curve
satisfying $I(\gamma)>\mu(\gamma)$?} 

\medskip
The authors do not know of any such examples.
If we suppose $I(\gamma)= \mu(\gamma)$, 
then 
\eqref{eq:hidden} yields the
inequality $I(\gamma)\le 2\#_\gamma$.

\medskip
\noindent
{Question 3.}
{\it Does 
$I(\gamma)\le 2\#_\gamma$ hold
for any generic closed curve in $\R^2$?}

\medskip
\noindent
{Question 4.}
{\it Is an arbitrary locally convex curve
$g$-tight?} 

\medskip
\noindent
{Question 5.}
{\it Find a criterion for $g$-tightness. 
For example, can one determine 
$\delta(\gamma)$ when 
$\gamma$ is of type $2_2$ or $1^b_2$?
}

\begin{acknowledgements}
 The authors thank Wayne Rossman  
 for careful reading of the first draft and for giving
valuable comments.
\end{acknowledgements}

\section{Tables of curves.}

\end{document}